\documentclass[12pt]{article}

\usepackage{amsfonts, amsthm, amsmath, amssymb}
\usepackage{breakurl}
\usepackage{color}
\usepackage{fancyvrb}
\usepackage{float}
\usepackage{graphicx}
\usepackage{hyperref}
\usepackage[nottoc]{tocbibind} 
\usepackage{url}

\allowdisplaybreaks

\hypersetup{
pdftitle={Two statements that are equivalent to a conjecture related to the distribution of prime numbers},
pdfauthor={Germ\'an Andr\'es Paz},
pdfsubject={},
pdfkeywords={Andrica's conjecture, Brocard's conjecture, Legendre's conjecture, Oppermann's conjecture, prime numbers, triangular numbers}
}

\makeatletter
\def\thmhead@plain#1#2#3{%
  \thmname{#1}\thmnumber{\@ifnotempty{#1}{ }\@upn{#2}}%
  \thmnote{ {\the\thm@notefont#3}}}
\let\thmhead\thmhead@plain
\makeatother

\newcommand{\A}{\operatorname{Set\;A}}

\newcommand{\B}{\operatorname{Set\;B}}

\newcommand{\xqedhere}[2]{%
  \rlap{\hbox to#1{\hfil\llap{\ensuremath{#2}}}}}

\theoremstyle{definition}

\newtheorem{Theorem}{Theorem}[section]
\newtheorem{Conjecture}[Theorem]{Conjecture}
\newtheorem{Corollary}[Theorem]{Corollary}
\newtheorem{Definition}[Theorem]{Definition}
\newtheorem{Lemma}[Theorem]{Lemma}
\newtheorem{Remark}[Theorem]{Remark}
\newtheorem{Statement}[Theorem]{Statement}

\newtheorem*{Andrica's conjecture}{Andrica's conjecture}
\newtheorem*{Brocard's conjecture}{Brocard's conjecture}
\newtheorem*{Legendre's conjecture}{Legendre's conjecture}
\newtheorem*{Oppermann's conjecture}{Oppermann's conjecture}

\begin{document}

\title{Two statements that are equivalent to a conjecture related to the distribution of prime numbers}
\author{Germ\'an Andr\'es Paz}
\date{June 19, 2014}

\maketitle

\centerline{\bf Abstract}

Let $n\in\mathbb{Z}^+$. In \cite{Paz} we ask the question whether any sequence of $n$ consecutive integers greater than $n^2$ and smaller than $(n+1)^2$ contains at least one prime number, and we show that this is actually the case for every $n\leq 1,193,806,023$. In addition, we prove that a positive answer to the previous question for all $n$ would imply Legendre's, Brocard's, Andrica's, and Oppermann's conjectures, as well as the assumption that for every $n$ there is always a prime number in the interval $[n,n+2\lfloor\sqrt{n}\rfloor-1]$.

Let $\pi[n+g(n),n+f(n)+g(n)]$ denote the amount of prime numbers in the interval $[n+g(n),n+f(n)+g(n)]$. Here we show that the conjecture described in \cite{Paz} is equivalent to the statement that
\begin{equation*}
\pi[n+g(n),n+f(n)+g(n)]\ge 1\text{, }\forall n\in\mathbb{Z}^+\text{,}
\end{equation*}
where
\begin{equation*}
f(n)=\left(\frac{n-\lfloor\sqrt{n}\rfloor^2-\lfloor\sqrt{n}\rfloor-\beta}{|n-\lfloor\sqrt{n}\rfloor^2-\lfloor\sqrt{n}\rfloor-\beta|}\right)(1-\lfloor\sqrt{n}\rfloor)\text{, }g(n)=\left\lfloor1-\sqrt{n}+\lfloor\sqrt{n}\rfloor\right\rfloor\text{,}
\end{equation*} 
and $\beta$ is any real number such that $1<\beta<2$. We also prove that the conjecture in question is equivalent to the statement that
\begin{equation*}
\pi[S_n,S_n+\lfloor\sqrt{S_n}\rfloor-1]\ge 1\text{, }\forall n\in\mathbb{Z}^+\text{,}
\end{equation*}
where
\begin{equation*}
S_n=n+\frac{1}{2}\left\lfloor\frac{\sqrt{8n+1}-1}{2}\right\rfloor^2-\frac{1}{2}\left\lfloor\frac{\sqrt{8n+1}-1}{2}\right\rfloor+1\text{.}
\end{equation*}
We use this last result in order to create plots of $h(n)=\pi[S_n,S_n+\lfloor\sqrt{S_n}\rfloor-1]$ for many values of $n$.



\section{Introduction}

The well-known Bertrand's postulate states that for every integer $n>3$ there always exists a prime number $p$ such that $n<p<2n-2$ (another formulation of this theorem is that for every $n>1$ there always exists a prime number $p$ such that $n<p<2n$). This statement, which had been conjectured by Joseph Bertrand in 1845, was first proved by P. L. Chebyshev in 1850.

After Bertrand's postulate was proved in 1850, many better results have been obtained by using both elementary and nonelementary methods. In 1930, Hoheisel \cite{Hoheisel} showed that there exists a constant $\theta<1$ such that 
\begin{equation}\label{equation0}
\pi(x+x^\theta)-\pi(x)\sim\frac{x^\theta}{\ln x}\text{,}
\end{equation}
where $\pi$ denotes the prime-counting function. In fact, Hoheisel showed that one may take $\theta=32999/33000$. This result has been subsequently improved, and in 2001 Baker, Harman, and Pintz \cite{Baker_Harman_Pintz} proved that in \eqref{equation0} the constant $\theta$ may be taken to be $0.525$. In other words, these authors showed that the interval $[x,x+x^{0.525}]$ contains at least one prime number for sufficiently large $x$. Besides, we can also say that if the Riemann hypothesis is true, then in \eqref{equation0} we can take $\theta=1/2+\varepsilon$ \cite{Maier}.

Although much progress has been made towards finding shorter and shorter intervals containing at least one prime number, there are still many unsolved problems in Number Theory regarding the existence of prime numbers in certain intervals. The following four conjectures are examples of such problems:

\begin{Legendre's conjecture}[\cite{Hardy_Wright, Landau, Legendre}]
For every positive integer $n$ there exists at least one prime number $p$ such that $n^2<p<(n+1)^2$.\hfill$\triangleleft$
\end{Legendre's conjecture}

Yitang Zhang \cite{Zhang} made a groundbreaking discovery regarding gaps between prime numbers. Zhang proved that
\begin{equation*}
\liminf_{n\to\infty}(p_{n+1}-p_n)<7\times 10^7\text{.}
\end{equation*}
This result has now been significantly improved (see, for example, \cite{Maynard} and \cite{Polymath8} for more details).

Legendre's conjecture has not been proved or disproved yet, but one thing we know thanks to Zhang's result is that there are infinitely many positive integers $n$ such that between $n^2$ and $(n+2)^2$ there are at least two prime numbers.

\begin{Brocard's conjecture}[\cite{Brocard}]
If $p_n$ and $p_{n+1}$ are consecutive prime numbers greater than 2, then between $(p_n)^2$ and $(p_{n+1})^2$ there are at least four prime numbers.\hfill$\triangleleft$
\end{Brocard's conjecture}

\begin{Andrica's conjecture}[\cite{Andrica1, Andrica2}]
The inequality $\sqrt{p_{n+1}}-\sqrt{p_n}<1$ holds for every pair of consecutive prime numbers $p_n$ and $p_{n+1}$.\hfill$\triangleleft$
\end{Andrica's conjecture}

Although Andrica's conjecture has not been proved (or disproved) yet, we know that Zhang's result \cite{Zhang} implies that there are infinitely many pairs of consecutive prime numbers $p_n$ and $p_{n+1}$ such that $\sqrt{p_{n+1}}-\sqrt{p_n}<1$.

\begin{Oppermann's conjecture}[\cite{Oppermann}]
For any integer $n>1$ there is a prime number in the interval $[n^2-n,n^2]$ and another prime in the interval $[n^2,n^2+n]$.\hfill$\triangleleft$
\end{Oppermann's conjecture}

In \cite{Paz} we consider the following conjecture regarding the distribution of prime numbers:

\begin{Conjecture}\label{Conjecture1}
If $n$ is any positive integer and we take $n$ consecutive integers
located between $n^2$ and $(n+1)^2$, then among those
$n$ integers there is at least one prime number. In other words, if $a_1$, $a_2$, $a_3$, $a_4$, \ldots{}, $a_n$ are $n$ consecutive integers
such that $n^2<a_1<a_2<a_3<a_4<...<a_n<(n + 1)^2$, then at least one of those $n$ integers is a prime number.\hfill$\triangleleft$
\end{Conjecture}

\begin{Remark}
Throughout this paper, whenever we say that a number $y$ is
\emph{between} a number $x$ and a number $z$, it means that
$x<y<z$, which means that $y$ is never equal to $x$ or
$z$.\hfill$\triangleleft$
\end{Remark}

While it is conjectured that \eqref{equation0} holds for all $\theta\in(0,1)$, it would be interesting to find a relation that explains why a sequence of $n$ consecutive integers greater than $n^2$ and smaller than $(n+1)^2$ cannot contain only composite numbers. In other words, it would be interesting to find a relation between the amount of consecutive integers in a sequence and the perfect squares between which the mentioned sequence is located.

In \cite{Paz} we prove that if Conjecture \ref{Conjecture1} is true, then Legendre's, Brocard's, Andrica's, and Oppermann's conjectures follow. In addition, we show that if the mentioned conjecture holds, then there is always a prime number in the interval $[n,n+2\lfloor\sqrt{n}\rfloor-1]$ for every positive integer $n$. Moreover, we prove that the conjecture in question is true at least for every positive integer $n\leq 1,193,806,023$. In this paper we provide two statements that are equivalent to Conjecture \ref{Conjecture1}.

Let us consider the following lemma:

\begin{Lemma}\label{Lemma0}
Suppose that $x$ is a certain positive integer. If Conjecture \ref{Conjecture1} is true, then the following statements are all true:

\begin{itemize}

\item If $x$ is a perfect square, then in the interval $[x+1,x+\lfloor\sqrt{x}\rfloor]$ there is at least one prime number.

\item If $x$ is not a perfect square and $\lfloor\sqrt{x}\rfloor^2<x\leq\lfloor\sqrt{x}\rfloor^2+\lfloor\sqrt{x}\rfloor+1<(\lfloor\sqrt{x}\rfloor+1)^2$, then in the interval $[x,x+\lfloor\sqrt{x}\rfloor-1]$ there is at least one prime number.

\item If $x$ is not a perfect square and $\lfloor\sqrt{x}\rfloor^2<\lfloor\sqrt{x}\rfloor^2+\lfloor\sqrt{x}\rfloor+2\leq x<(\lfloor\sqrt{x}\rfloor+1)^2$, then in the interval $[x-\lfloor\sqrt{x}\rfloor+1,x]$ there is at least one prime number.\qed

\end{itemize}
\end{Lemma}

\begin{Remark}
The symbol $\lfloor{}\rfloor$ represents the \emph{floor function}. The floor function of a given number is the largest integer that is not greater than that number. For example, $\lfloor{3.5}\rfloor=3$.\hfill$\triangleleft$
\end{Remark}

Unless we know the value of $x$ in advance, we cannot know which of the three intervals $[x+1,x+\lfloor\sqrt{x}\rfloor]$, $[x,x+\lfloor\sqrt{x}\rfloor-1]$, and $[x-\lfloor\sqrt{x}\rfloor+1,x]$ must contain a prime if Conjecture \ref{Conjecture1} is true. For this reason, we want to find a unique interval that can be applied to \emph{any} positive integer.

\section{An interval with a nonperfect square as a parameter}

Suppose that $a$ is a positive integer and that $a$ is not a perfect square. This means that $a$ is located between two consecutive perfect squares $\lfloor\sqrt{a}\rfloor^2$ and $(\lfloor\sqrt{a}\rfloor+1)^2$. In other words, we have
\begin{equation*}
\lfloor\sqrt{a}\rfloor^2<a<(\lfloor\sqrt{a}\rfloor+1)^2\text{.}
\end{equation*}
It is very easy to prove that the amount of integers that are greater than $\lfloor\sqrt{a}\rfloor^2$ and smaller than $(\lfloor\sqrt{a}\rfloor+1)^2$ is equal to $2\lfloor\sqrt{a}\rfloor$.

Now we need to define some concepts:

\begin{Definition}\label{Definition1}
The set made up of all integers contained in the interval $[\lfloor\sqrt{a}\rfloor^2+1,\lfloor\sqrt{a}\rfloor^2+\lfloor\sqrt{a}\rfloor+1]$ will be denoted by $\A(\lfloor\sqrt{a}\rfloor^2)$.\hfill$\triangleleft$
\end{Definition}

\begin{Definition}
The set made up of all integers that belong to the interval $[\lfloor\sqrt{a}\rfloor^2+\lfloor\sqrt{a}\rfloor+2,(\lfloor\sqrt{a}\rfloor+1)^2-1]$ will be denoted by $\B(\lfloor\sqrt{a}\rfloor^2)$.\hfill$\triangleleft$
\end{Definition}

The largest integer that belongs to $\A(\lfloor\sqrt{a}\rfloor^2)$ is $\lfloor\sqrt{a}\rfloor^2+\lfloor\sqrt{a}\rfloor+1$, whereas the smallest integer that belongs to $\B(\lfloor\sqrt{a}\rfloor^2)$ is $\lfloor\sqrt{a}\rfloor^2+\lfloor\sqrt{a}\rfloor+2$. Now, if $\alpha$ is any real number such that $0<\alpha<1$, we have
\begin{equation*}
\lfloor\sqrt{a}\rfloor^2+\lfloor\sqrt{a}\rfloor+1<\lfloor\sqrt{a}\rfloor^2+\lfloor\sqrt{a}\rfloor+1+\alpha<\lfloor\sqrt{a}\rfloor^2+\lfloor\sqrt{a}\rfloor+2\text{.}
\end{equation*}
We can also say that if $\beta$ is any real number such that $1<\beta<2$ (that is to say, if $\beta=1+\alpha$), then we have
\begin{equation}\label{equation1}
\lfloor\sqrt{a}\rfloor^2+\lfloor\sqrt{a}\rfloor+1<\lfloor\sqrt{a}\rfloor^2+\lfloor\sqrt{a}\rfloor+\beta<\lfloor\sqrt{a}\rfloor^2+\lfloor\sqrt{a}\rfloor+2\text{.}
\end{equation}
Taking \eqref{equation1} into account, we can state the following two lemmas:

\begin{Lemma}\label{Lemma1}
If $a\in\A(\lfloor\sqrt{a}\rfloor^2)$ (that is to say, if $a$ belongs to $\A(\lfloor\sqrt{a}\rfloor^2)$), then $a-(\lfloor\sqrt{a}\rfloor^2+\lfloor\sqrt{a}\rfloor+\beta)=a-\lfloor\sqrt{a}\rfloor^2-\lfloor\sqrt{a}\rfloor-\beta<0$.\qed
\end{Lemma}

\begin{Lemma}\label{Lemma2}
If $a\in\B(\lfloor\sqrt{a}\rfloor^2)$, then $a-(\lfloor\sqrt{a}\rfloor^2+\lfloor\sqrt{a}\rfloor+\beta)=a-\lfloor\sqrt{a}\rfloor^2-\lfloor\sqrt{a}\rfloor-\beta>0$.\qed
\end{Lemma}

Lemmas \ref{Lemma1} and \ref{Lemma2} lead to the following two lemmas:

\begin{Lemma}\label{Lemma3}
If $a\in\A(\lfloor\sqrt{a}\rfloor^2)$, then
\begin{equation*}
\pushQED{\qed} 
\frac{a-\lfloor\sqrt{a}\rfloor^2-\lfloor\sqrt{a}\rfloor-\beta}{|a-\lfloor\sqrt{a}\rfloor^2-\lfloor\sqrt{a}\rfloor-\beta|}=-1\text{.}\qedhere
\popQED
\end{equation*}
\end{Lemma}

\begin{Lemma}\label{Lemma4}
If $a\in\B(\lfloor\sqrt{a}\rfloor^2)$, then
\begin{equation*}
\pushQED{\qed}
\frac{a-\lfloor\sqrt{a}\rfloor^2-\lfloor\sqrt{a}\rfloor-\beta}{|a-\lfloor\sqrt{a}\rfloor^2-\lfloor\sqrt{a}\rfloor-\beta|}=1\text{.}\qedhere
\popQED
\end{equation*}
\end{Lemma}

\begin{Remark}
Throughout this paper, we use the symbol $||$ to denote the \emph{absolute value} of a number or expression. For example, the absolute value of $x$ is denoted by $|x|$.\hfill$\triangleleft$
\end{Remark}

If we consider Lemma \ref{Lemma0}, we can state the following two lemmas:

\begin{Lemma}\label{Lemma5}
If $a\in\A(\lfloor\sqrt{a}\rfloor^2)$ and Conjecture \ref{Conjecture1} is true, then in the interval $[a,a+\lfloor\sqrt{a}\rfloor-1]$ there is at least one prime number.\qed
\end{Lemma}

\begin{Lemma}\label{Lemma6}
If $a\in\B(\lfloor\sqrt{a}\rfloor^2)$ and Conjecture \ref{Conjecture1} is true, then in the interval $[a-\lfloor\sqrt{a}\rfloor+1,a]$ there is at least one prime number.\qed
\end{Lemma}

Combining Lemmas \ref{Lemma3}, \ref{Lemma4}, \ref{Lemma5}, and \ref{Lemma6}, we obtain the following result, which we will express as a theorem:

\begin{Theorem}\label{Theorem1}
Suppose that $a\in\mathbb{Z}^+$ and that $a$ is not a perfect square. If Conjecture \ref{Conjecture1} is true and $\beta$ is any real number such that $1<\beta<2$, then the interval
\begin{equation}\label{equation2}
\left[a,a+\left(\frac{a-\lfloor\sqrt{a}\rfloor^2-\lfloor\sqrt{a}\rfloor-\beta}{|a-\lfloor\sqrt{a}\rfloor^2-\lfloor\sqrt{a}\rfloor-\beta|}\right)(1-\lfloor\sqrt{a}\rfloor)\right]
\end{equation}
contains at least one prime number.\qed
\end{Theorem}

\section{An interval with a perfect square as a parameter}

Suppose that $b$ is a positive integer and that $b$ is a perfect square. If Conjecture \ref{Conjecture1} is true, then in the interval
\begin{equation*}
[b,b+\sqrt{b}]
\end{equation*}
there is at least one prime number. Since $b$ is a perfect square, it follows that
\begin{equation*}
\sqrt{b}=\lfloor\sqrt{b}\rfloor\text{.}
\end{equation*}
Therefore, we can also say that if Conjecture \ref{Conjecture1} is true, then there is always a prime number in the interval
\begin{equation*}
[b,b+\lfloor\sqrt{b}\rfloor]\text{.}
\end{equation*}

The interval \eqref{equation2} can be applied to every positive integer $a$ such that $a$ is not a perfect square, but it cannot be applied to $b$, since $b$ is a perfect square. Note that if Conjecture \ref{Conjecture1} is true, then the interval $[b,b+\lfloor\sqrt{b}\rfloor]$ contains at least one prime, but since $b$ is a perfect square we have
\begin{equation*}
b+\left(\frac{b-\lfloor\sqrt{b}\rfloor^2-\lfloor\sqrt{b}\rfloor-\beta}{|b-\lfloor\sqrt{b}\rfloor^2-\lfloor\sqrt{b}\rfloor-\beta|}\right)(1-\lfloor\sqrt{b}\rfloor)=b+\lfloor\sqrt{b}\rfloor-1\text{,}
\end{equation*}
which means that the interval
\begin{equation}\label{equation3}
\left[b,b+\left(\frac{b-\lfloor\sqrt{b}\rfloor^2-\lfloor\sqrt{b}\rfloor-\beta}{|b-\lfloor\sqrt{b}\rfloor^2-\lfloor\sqrt{b}\rfloor-\beta|}\right)(1-\lfloor\sqrt{b}\rfloor)\right]
\end{equation}
does not contain any prime in the case where $b=1$. However, we can say that if $b$ is \emph{any} perfect square and Conjecture \ref{Conjecture1} is true, then in the interval
\begin{equation*}
\left[b,b+\left(\frac{b-\lfloor\sqrt{b}\rfloor^2-\lfloor\sqrt{b}\rfloor-\beta}{|b-\lfloor\sqrt{b}\rfloor^2-\lfloor\sqrt{b}\rfloor-\beta|}\right)(1-\lfloor\sqrt{b}\rfloor)+1\right]
\end{equation*}
there is at least one prime number. We can see that we are adding 1 to the upper endpoint of the interval \eqref{equation3} in case $b$ is set to 1. In fact, we can also add 1 to the lower endpoint of \eqref{equation3}, since no square of a positive integer is a prime number. 

We are now ready to state the following theorem:

\begin{Theorem}\label{Theorem2}
Suppose that $b\in\mathbb{Z}^+$ and that $b$ is a perfect square. If Conjecture \ref{Conjecture1} is true and $\beta$ is any real number such that $1<\beta<2$, then the interval
\begin{equation*}
\left[b+1,b+\left(\frac{b-\lfloor\sqrt{b}\rfloor^2-\lfloor\sqrt{b}\rfloor-\beta}{|b-\lfloor\sqrt{b}\rfloor^2-\lfloor\sqrt{b}\rfloor-\beta|}\right)(1-\lfloor\sqrt{b}\rfloor)+1\right]
\end{equation*}
contains at least one prime number.\qed
\end{Theorem}

\section{A unique interval that can have any positive integer as a parameter}

We start this section by stating the following lemma:

\begin{Lemma}\label{Lemma7}
Suppose that $x$ is any positive integer and consider the function $g(x)=\left\lfloor1-\sqrt{x}+\lfloor\sqrt{x}\rfloor\right\rfloor$. The function $g(x)$ has the following property:
\begin{equation*}
g(x) = 
\begin{cases}
0\text{,} &\mbox{if } x \text{ is not a perfect square;} \\ 
1\text{,} &\mbox{if } x \text{ is a perfect square.} \xqedhere{112.5pt}{\qed} 
\end{cases}
\end{equation*}
\end{Lemma}

If we consider Theorems \ref{Theorem1} and \ref{Theorem2} and Lemma \ref{Lemma7}, we can state the following theorem:

\begin{Theorem}\label{Theorem3}
Suppose that $n$ is any positive integer and $\beta$ is any real number such that $1<\beta<2$. If Conjecture \ref{Conjecture1} is true, then the interval
\begin{multline*}
\bigg[n+\left\lfloor1-\sqrt{n}+\lfloor\sqrt{n}\rfloor\right\rfloor,\\n+\left(\frac{n-\lfloor\sqrt{n}\rfloor^2-\lfloor\sqrt{n}\rfloor-\beta}{|n-\lfloor\sqrt{n}\rfloor^2-\lfloor\sqrt{n}\rfloor-\beta|}\right)(1-\lfloor\sqrt{n}\rfloor)+\left\lfloor1-\sqrt{n}+\lfloor\sqrt{n}\rfloor\right\rfloor\bigg]
\end{multline*}
contains at least one prime number.\qed
\end{Theorem}

The number
\begin{equation*}
n+\left(\frac{n-\lfloor\sqrt{n}\rfloor^2-\lfloor\sqrt{n}\rfloor-\beta}{|n-\lfloor\sqrt{n}\rfloor^2-\lfloor\sqrt{n}\rfloor-\beta|}\right)(1-\lfloor\sqrt{n}\rfloor)+\left\lfloor1-\sqrt{n}+\lfloor\sqrt{n}\rfloor\right\rfloor
\end{equation*}
is never prime when $n$ is a perfect square greater than 1. As a consequence, the following corollary is deduced from Theorem \ref{Theorem3}:

\begin{Corollary}
Suppose that $n$ is any positive integer greater than 1 and $\beta$ is any real number such that $1<\beta<2$. If Conjecture \ref{Conjecture1} is true, then the interval
\begin{equation*}
\bigg[n+\left\lfloor1-\sqrt{n}+\lfloor\sqrt{n}\rfloor\right\rfloor,n+\left(\frac{n-\lfloor\sqrt{n}\rfloor^2-\lfloor\sqrt{n}\rfloor-\beta}{|n-\lfloor\sqrt{n}\rfloor^2-\lfloor\sqrt{n}\rfloor-\beta|}\right)(1-\lfloor\sqrt{n}\rfloor)\bigg]
\end{equation*}
contains at least one prime number.\qed
\end{Corollary}

\section{Another interval that can have any positive integer as a parameter}

We start this section with the following definition:

\begin{Definition}
Let us look at the following sequence of numbers:

\begin{center}
$\mathbf{1^2}$\quad \textcolor{blue}{2}\quad \textcolor{blue}{3}\quad $\mathbf{2^2}$\quad \textcolor{blue}{5}\quad \textcolor{blue}{6}\quad \textcolor{blue}{7}\quad 8\quad $\mathbf{3^2}$\quad \textcolor{blue}{10}\quad \textcolor{blue}{11}\quad \textcolor{blue}{12}\quad \textcolor{blue}{13}\quad 14\quad 15\quad $\mathbf{4^2}$
\end{center}

\begin{center}
\textcolor{blue}{17}\quad \textcolor{blue}{18}\quad \textcolor{blue}{19}\quad \textcolor{blue}{20}\quad \textcolor{blue}{21}\quad 22\quad 23\quad 24\quad $\mathbf{5^2}$\quad \textcolor{blue}{26}\quad \textcolor{blue}{27}\quad \textcolor{blue}{28}\quad \textcolor{blue}{29}\quad \textcolor{blue}{30}\quad \textcolor{blue}{31}\quad $\cdots$
\end{center}
The numbers in blue are the integers in $\A(1^2)$ (see Definition \ref{Definition1}), the numbers in $\A(2^2)$, the numbers in $\A(3^2)$, etc. These numbers form a sequence which we will denote by $S$. In other words,
\begin{equation*}
S=(2,3,5,6,7,10,11,12,13,17,18,19,20,21,26,27,28,29,30,31,\dots{})\text{.}\xqedhere{13.5pt}{\triangleleft}
\end{equation*}
\end{Definition}

In this section we will find an explicit formula that allows us to calculate $S_n$, that is to say, the $n$th term of sequence $S$. It is obvious that Conjecture \ref{Conjecture1} is equivalent to the statement that for every positive integer $n$ there is always a prime number in the interval $[S_n,S_n+\lfloor\sqrt{S_n}\rfloor-1]$ (see Lemma \ref{Lemma5}).  

Now, let us take the following definitions into account:

\begin{Definition}
A positive integer is said to be a \emph{triangular number} if it is a number of the form
\begin{equation*}
\sum_{k=1}^{m}k\text{,}
\end{equation*}
where $k$ and $m$ are positive integers and $m\ge k$. The sequence of triangular numbers will be denoted by $T$.\hfill$\triangleleft$
\end{Definition}

\begin{Definition}
The sequence of positive integers that are equal to a perfect square plus 1 will be denoted by $U$.\hfill$\triangleleft$
\end{Definition}

In order to find the $n$th triangular number, we can use the formula
\begin{equation}\label{equation4}
T_n=(n^2+n)/2\text{;}
\end{equation}
in order to find the $n$th term of sequence $U$, we can use the formula
\begin{equation}\label{equation5}
U_n=n^2+1\text{.}
\end{equation}
Combining \eqref{equation4} and \eqref{equation5}, we get
\begin{equation}\label{equation6}
U_n-T_n=(n^2-n)/2+1\text{.}
\end{equation}

Now we will arrange the sequence of positive integers and sequence $S$ in the following way:
\vspace{-12pt}
\begin{multline*}
\\
\text{Group 1}
\begin{cases}
\textcolor{white}{0}\textcolor{red}{1}\rightarrow \textcolor{white}{0}\textcolor{green}{2}\\ 
\textcolor{white}{0}2\rightarrow \textcolor{white}{0}3
\end{cases}\\
\text{Group 2}
\begin{cases}
\textcolor{white}{0}\textcolor{red}{3}\rightarrow \textcolor{white}{0}\textcolor{green}{5}\\ 
\textcolor{white}{0}4\rightarrow \textcolor{white}{0}6\\
\textcolor{white}{0}5\rightarrow \textcolor{white}{0}7
\end{cases}\\
\text{Group 3}
\begin{cases}
\textcolor{white}{0}\textcolor{red}{6}\rightarrow \textcolor{green}{10}\\ 
\textcolor{white}{0}7\rightarrow 11\\
\textcolor{white}{0}8\rightarrow 12\\
\textcolor{white}{0}9\rightarrow 13
\end{cases}\\
\text{Group 4}
\begin{cases}
\textcolor{red}{10}\rightarrow \textcolor{green}{17}\\ 
11\rightarrow 18\\
12\rightarrow 19\\
13\rightarrow 20\\
14\rightarrow 21
\end{cases}\\
\end{multline*}

\vspace{-27pt}
\hspace{137.5pt}\vdots\hspace{41pt}\vdots\hspace{11pt}\vdots\hspace{12.5pt}\vdots
\vspace{12pt}

\noindent In our graphic, the numbers in red are triangular numbers, whereas the numbers in green are integers that are equal to a perfect square plus 1. Now, let us take the following statements into account:

\begin{Statement}
The first column of numbers (from left to right) is the sequence of positive integers, whereas the second column is sequence $S$.\hfill$\triangleleft$
\end{Statement}

\begin{Statement}
Group 1 contains two rows of numbers, Group 2 contains three, Group 3 contains four, etc. In general, if $x$ is any positive integer, then Group $x$ will contain $x+1$ rows of numbers.\hfill$\triangleleft$
\end{Statement}

\begin{Statement}\label{Statement1}
In Group 1 the first line contains $T_1$ (in red) and $U_1$ (in green). In general, in Group $x$ the first line contains $T_x$ and $U_x$.\hfill$\triangleleft$
\end{Statement}

\begin{Statement}\label{Statement2}
In Group 1 the difference between two numbers that are in the same row is equal to $U_1-T_1$. In general, in Group $x$ the difference between two numbers in the same row is equal to $U_x-T_x$.\hfill$\triangleleft$
\end{Statement}

\begin{Statement}
The expression $1\rightarrow 2$ means that the first term of sequence $S$ is 1. In other words, the expression $1\rightarrow 2$ means that $S_1=2$.\hfill$\triangleleft$ 
\end{Statement}

Suppose we choose any positive integer. Is there a general formula that can tell us in which Group we will find that integer in the first column? The answer is yes. In order to find such formula, we need to take into account that if $n$ is a positive integer, then $n$ appears in the first column of Group $x$ if $T_x\leq n<T_{x+1}$ and vice versa.

We know that
\begin{equation*}
\frac{\sqrt{8n+1}-1}{2}
\end{equation*}
is an integer only if $n$ is a triangular number (see \cite{Triangular_Numbers}). Therefore, the fact that
\begin{equation*}
T_x\leq n<T_{x+1}
\end{equation*}
implies that
\begin{equation*}
x\leq\frac{\sqrt{8n+1}-1}{2}<x+1
\end{equation*}
and vice versa. This means that
\begin{equation}\label{equation7}
x=\left\lfloor\frac{\sqrt{8n+1}-1}{2}\right\rfloor\text{.}
\end{equation}
Combining \eqref{equation6}, \eqref{equation7}, and Statements \ref{Statement1} and \ref{Statement2}, we conclude that
\begin{align*}
S_n&=n+\left(\left\lfloor\frac{\sqrt{8n+1}-1}{2}\right\rfloor^2-\left\lfloor\frac{\sqrt{8n+1}-1}{2}\right\rfloor\right)/2+1\\
&=n+\frac{1}{2}\left\lfloor\frac{\sqrt{8n+1}-1}{2}\right\rfloor^2-\frac{1}{2}\left\lfloor\frac{\sqrt{8n+1}-1}{2}\right\rfloor+1\text{.}
\end{align*}
We will state this result as a lemma:

\begin{Lemma}\label{Lemma8}
If $n$ is any positive integer, then the $n$th term of sequence $S$ can be found using the following explicit formula:
\begin{equation*}
\pushQED{\qed}
S_n=n+\frac{1}{2}\left\lfloor\frac{\sqrt{8n+1}-1}{2}\right\rfloor^2-\frac{1}{2}\left\lfloor\frac{\sqrt{8n+1}-1}{2}\right\rfloor+1\text{.}\qedhere
\popQED
\end{equation*}
\end{Lemma}

As we said before, Conjecture \ref{Conjecture1} is equivalent to the statement that for every positive integer $n$ there is always a prime number in the interval $[S_n,S_n+\lfloor\sqrt{S_n}\rfloor-1]$. So, if we take Lemmas \ref{Lemma5} and \ref{Lemma8} into account, we can state the following theorem:

\begin{Theorem}\label{Theorem4}
Conjecture \ref{Conjecture1} is equivalent to the statement that
\begin{equation*}
\pi[S_n,S_n+\lfloor\sqrt{S_n}\rfloor-1]\ge 1\text{, }\forall n\in\mathbb{Z}^+\text{,}
\end{equation*}
where
\begin{equation*}
\pushQED{\qed}
S_n=n+\frac{1}{2}\left\lfloor\frac{\sqrt{8n+1}-1}{2}\right\rfloor^2-\frac{1}{2}\left\lfloor\frac{\sqrt{8n+1}-1}{2}\right\rfloor+1\text{.}\qedhere
\popQED
\end{equation*}
\end{Theorem}

Additionally, we will state the following lemma:

\begin{Lemma}
If Conjecture \ref{Conjecture1} is true, then for every $n$ there always exist two prime numbers $p$ and $q$ such that $2T_n<p<q<2T_{n+1}$.
\end{Lemma}

Zhang's theorem on bounded gaps between primes implies that there are infinitely many triangular numbers $T_n$ such that between $2T_n$ and $2T_{n+2}$ there are at least two prime numbers.

\section{Conclusion}

We have proved the following theorems:

\vspace{\topsep}

\noindent{\bf Theorem \ref{Theorem3}.} 
Let $\pi[n+g(n),n+f(n)+g(n)]$ denote the amount of prime numbers in the interval $[n+g(n),n+f(n)+g(n)]$. Conjecture \ref{Conjecture1} is equivalent to the statement that
\begin{equation*}
\pi[n+g(n),n+f(n)+g(n)]\ge1\text{, }\forall n\in\mathbb{Z}^+\text{,}
\end{equation*}
where
\begin{equation*}
f(n)=\left(\frac{n-\lfloor\sqrt{n}\rfloor^2-\lfloor\sqrt{n}\rfloor-\beta}{|n-\lfloor\sqrt{n}\rfloor^2-\lfloor\sqrt{n}\rfloor-\beta|}\right)(1-\lfloor\sqrt{n}\rfloor)\text{, }g(n)=\left\lfloor1-\sqrt{n}+\lfloor\sqrt{n}\rfloor\right\rfloor\text{,}
\end{equation*} 
and $\beta$ is any real number such that $1<\beta<2$.\qed

\vspace{\topsep}

\noindent{\bf Theorem \ref{Theorem4}.}
Conjecture \ref{Conjecture1} is equivalent to the statement that
\begin{equation*}
\pi[S_n,S_n+\lfloor\sqrt{S_n}\rfloor-1]\ge 1\text{, }\forall n\in\mathbb{Z}^+\text{,}
\end{equation*}
where
\begin{equation*}
\pushQED{\qed}
S_n=n+\frac{1}{2}\left\lfloor\frac{\sqrt{8n+1}-1}{2}\right\rfloor^2-\frac{1}{2}\left\lfloor\frac{\sqrt{8n+1}-1}{2}\right\rfloor+1\text{.}\qedhere
\popQED
\end{equation*}

\appendix

\section{Plots}

We define the function $h(n)$ in the following way:
\begin{equation*}
h(n)=\pi[S_n,S_n+\lfloor\sqrt{S_n}\rfloor-1]\text{.}
\end{equation*}
Now we will use Wolfram \emph{Mathematica} (version 9) in order to create plots of $h(n)$ for many values of $n$. In order to do so, we will take the following lemma into account:

\begin{Lemma}
Let $x,y\in\mathbb{Z}^+$ and let $x\leq y$. If $\pi(x)$ is the amount of prime numbers that are less than or equal to $x$, $\pi(y)$ the amount of prime numbers that are less than or equal to $y$, and $\pi[x,y]$ the amount of prime numbers in the interval $[x,y]$, then
\begin{equation*}
\pi[x,y] =
\begin{cases}
\pi(y)-\pi(x)+1\text{,} &\mbox{if } x \text{ is a prime number;} \\ 
\pi(y)-\pi(x)\text{,} &\mbox{otherwise.} \xqedhere{126.5pt}{\qed} 
\end{cases}
\end{equation*}
\end{Lemma}

In order to plot $h(n)$ for all $n\leq 100$, we can enter the following code in \emph{Mathematica}:
\begin{center}
\begin{minipage}{12.9cm}
\begin{Verbatim}[fontsize=\small, frame=single]
NumPrimes[m_Integer, n_Integer] := If[PrimeQ[m],
   PrimePi[n] - PrimePi[m] + 1,
   PrimePi[n] - PrimePi[m]
];

S[n_] := n + 1/2 \[LeftFloor](Sqrt[8 n + 1] - 1)/2\[RightFloor]
^2 - 1/2 \[LeftFloor](Sqrt[8 n + 1] - 1)/2\[RightFloor] + 1;

ListPlot[Table[NumPrimes[S[n], S[n] + \[LeftFloor]Sqrt[S[n]]
\[RightFloor] - 1], {n, 1, 100}], 
PlotRange -> {{0, 101}, {0, 4.5}}, Filling -> Axis]
\end{Verbatim}
\end{minipage}
\end{center}
\noindent If we copy the code directly from this PDF file, we may need to remove the extra space between \texttt{[RightFloor]} and \texttt{\^} in order for us not to get any error messages.

The code above gives us the following result:
\begin{figure}[H]
\begin{center}
\frame{\includegraphics[scale=0.4]{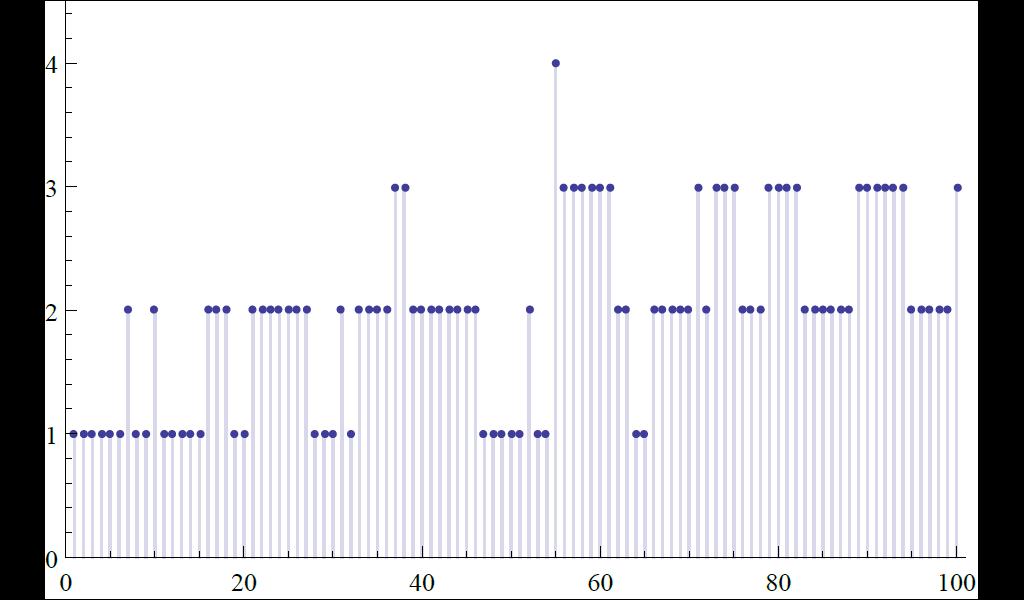}}
\caption{The function $h(n)=\pi[S_n,S_n+\lfloor\sqrt{S_n}\rfloor-1]$ for $n\leq 100$.}
\end{center}
\end{figure}
\noindent We can modify the code in order to create plots of $h(n)$ for more values of $n$. If we remove \texttt{, PlotRange -> \{\{0, 101\}, \{0, 4.5\}\}, Filling -> Axis} and we replace \texttt{\{n, 1, 100\}} with \texttt{\{n, 1, 10\;000\}} and with \texttt{\{n, 1, 500\;000\}}, we obtain the following plots:
\begin{figure}[H]
\begin{center}
\frame{\includegraphics[scale=0.4]{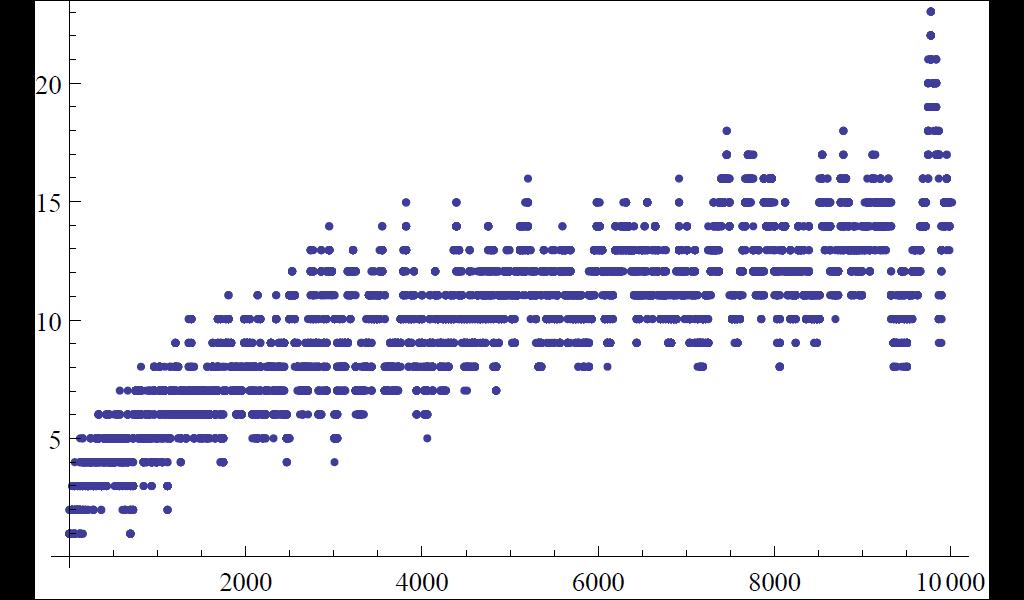}}
\caption{The function $h(n)=\pi[S_n,S_n+\lfloor\sqrt{S_n}\rfloor-1]$ for $n\leq 10,000$.}
\end{center}
\end{figure}
\begin{figure}[H]
\begin{center}
\frame{\includegraphics[scale=0.4]{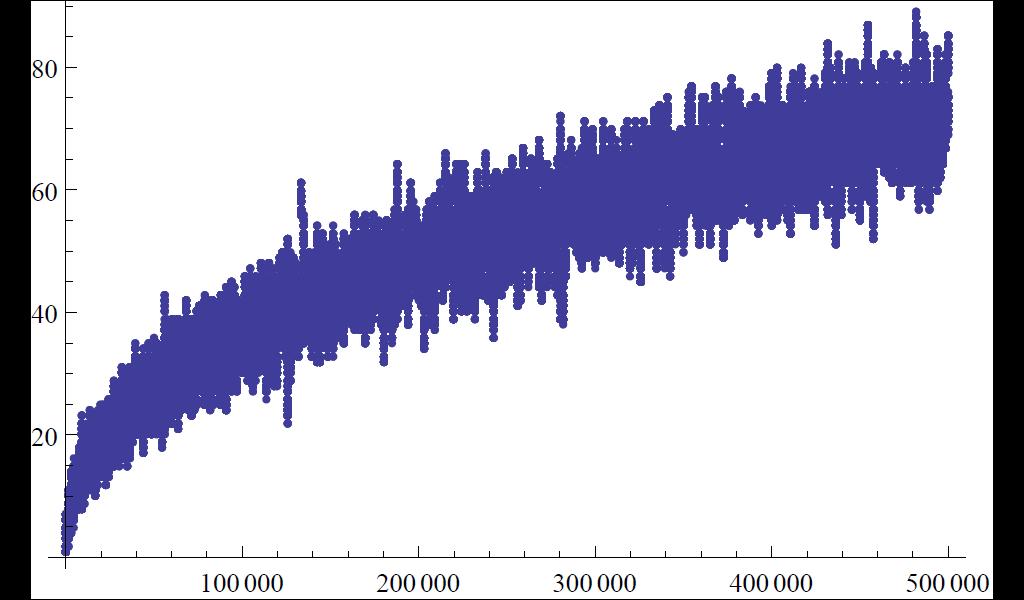}}
\caption{The function $h(n)=\pi[S_n,S_n+\lfloor\sqrt{S_n}\rfloor-1]$ for $n\leq 500,000$.}
\end{center}
\end{figure}


\noindent\emph{Instituto de Educaci\'on Superior N$^\circ$28 Olga Cossettini, (2000) Rosario, Santa Fe, Argentina\\germanpaz\_ar@hotmail.com}

\end{document}